# Ratio Estimators in Simple Random Sampling Using Information on Auxiliary Attribute

Rajesh Singh
Department of Statistics, Banaras Hindu University
Varanasi (U.P.), India
Email: rsinghstat@yahoo.co.in

Pankaj Chauhan
School of Statistics, DAVV, Indore (M.P.), India

Nirmala Sawan
School of Statistics, DAVV, Indore (M.P.), India

Florentin Smarandache
Department of Mathematics, University of New Mexico, Gallup, USA
Email: smarand@unm.edu

## Abstract

Some ratio estimators for estimating the population mean of the variable under study, which make use of information regarding the population proportion possessing certain attribute, are proposed. Under simple random sampling without replacement (SRSWOR) scheme, the expressions of bias and mean-squared error (MSE) up to the first order of approximation are derived. The results obtained have been illustrated numerically by taking some empirical population considered in the literature.

## 1. Introduction

The use of auxiliary information can increase the precision of an estimator when study variable y is highly correlated with auxiliary variable x. There exist situations when information is available in the form of attribute $\phi$, which is highly correlated with y. For example

a) Sex and height of the persons,
b) Amount of milk produced and a particular breed of the cow,
c) Amount of yield of wheat crop and a particular variety of wheat etc. (see Jhajj et. al. [1]).

Consider a sample of size n drawn by SRSWOR from a population of size N. Let $y_i$ and $\phi_i$ denote the observations on variable y and $\phi$ respectively for $i^{th}$ unit $(i = 1,2,....N)$. Suppose there is a complete dichotomy in the population with respect to the presence or absence of an attribute, say $\phi$, and it is assumed that attribute $\phi$ takes only the two values 0 and 1 according as

$$\phi_i = 1, \text{ if ith unit of the population possesses attribute } \phi$$
$$= 0, \text{ otherwise.}$$

Let $A=\sum_{i=1}^{N}\phi_i$ and $a=\sum_{i=1}^{n}\phi_i$ denote the total number of units in the population and sample respectively possessing attribute $\phi$. Let $P=\frac{A}{N}$ and $p=\frac{a}{n}$ denote the proportion of units in the population and sample respectively possessing attribute $\phi$.

Taking into consideration the point biserial correlation between a variable and an attribute, Naik and Gupta [2] defined ratio estimator of population mean when the prior information of population proportion of units, possessing the same attribute is available, as follows:

$$t_{NG}=\bar{y}\left(\frac{P}{p}\right) \tag{1.1}$$

Here $\bar{y}$ is the sample mean of variable of interest. The MSE of $t_{NG}$ up to the first order of approximation is –

$$MSE(t_{NG})=\left(\frac{1-f}{n}\right)\left[S_y^2+R_1^2S_\phi^2-2R_1S_{y\phi}\right] \tag{1.2}$$

Where $f=\frac{n}{N}$, $R_1=\frac{\bar{Y}}{P}$, $S_y^2=\frac{1}{N-1}\sum_{i=1}^{N}\left(y_i-\bar{Y}\right)^2$, $S_\phi^2=\frac{1}{N-1}\sum_{i=1}^{N}\left(\phi_i-P\right)^2$, $S_{y\phi}=\frac{1}{N-1}\sum_{i=1}^{N}\left(\phi_i-P\right)\left(y_i-\bar{Y}\right)$.

In the present paper, some ratio estimators for estimating the population mean of the variable under study, which make use of information regarding the population proportion possessing certain attribute, are proposed. The expressions of bias and MSE have been obtained. The numerical illustrations have also been done by taking some empirical populations considered in the literature.

## 2. The suggested estimator

Following Ray and Singh [3], we propose the following estimator –

$$t_1=\frac{\bar{y}+b_\phi(P-p)}{p}P=R^*P \tag{2.1}$$

Where $b_\phi=\frac{s_{y\phi}}{s_\phi^2}$, $R^*=\frac{\bar{y}+b_\phi(P-p)}{p}$, $s_\phi^2=\left(\frac{1}{n-1}\right)\sum_{i=1}^{n}\left(\phi_i-p\right)^2$ and $s_{y\phi}=\left(\frac{1}{n-1}\right)\sum_{i=1}^{n}\left(\phi_i-p\right)\left(y_i-\bar{Y}\right)$.

**Remark 1**: When we put $b_\phi=0$ in (2.1) the proposed estimator turns to the Naik and Gupta [2] ratio estimator $t_{NG}$ given in (1.1).

MSE of this estimator can be found by using Taylor series expansion given by –

$$f(p,\bar{y}) \cong f(P,\bar{y}) + \frac{\partial f(c,d)}{\partial c}\Big|_{P,\bar{Y}}(p-P) + \frac{\partial f(c,d)}{\partial c}\Big|_{P,\bar{Y}}(\bar{y}-\bar{Y}) \qquad (2.2)$$

Where $f(p,\bar{y}) = R^*$ and $f(P,\bar{Y}) = R_1$.

Expression (2.2) can be applied to the proposed estimator in order to obtain MSE equation as follows:

$$R^* - R_1 \cong \frac{\partial\left(\left(\bar{y}+b_\phi(P-p)\right)\right)/p}{\partial p}\Bigg|_{P,\bar{Y}}\left(p-P\right) + \frac{\partial\left(\left(\bar{y}+b_\phi(P-p)\right)\right)/p}{\partial \bar{y}}\Bigg|_{P,\bar{Y}}\left(\bar{y}-\bar{Y}\right)$$

$$\cong -\left(\frac{\bar{y}}{p^2}+\frac{b_\phi P}{p^2}\right)\Bigg|_{P,\bar{Y}}(p-P) + \frac{1}{p}\Bigg|_{P,\bar{Y}}(\bar{y}-\bar{Y})$$

$$E(R^*-R_1)^2 \cong \frac{\left(\bar{Y}+B_\phi P\right)^2}{P^4}V(p) - \frac{2(\bar{Y}+B_\phi P)}{P^3}Cov(p,\bar{y}) + \frac{1}{P^2}V(\bar{y})$$

$$\cong \frac{1}{P^2}\left\{\frac{(\bar{Y}+B_\phi P)^2}{P^2}V(p) - \frac{2(\bar{Y}+B_\phi P)}{P}Cov(p,\bar{y}) + V(\bar{y})\right\} \qquad (2.3)$$

Where $B_\phi = \frac{S_{\phi y}}{S_\phi^2} = \frac{\rho_{pb}S_y}{S_\phi}$.

$\rho_{pb} = \frac{S_{y\phi}}{S_y S_\phi}$, is the point biserial correlation coefficient.

Now,

$$MSE(t_1) = P^2E(R_1 - R_\phi)^2$$

$$\cong \frac{\left(\bar{Y}+B_\phi P\right)^2}{P^2}V\left(p\right) - \frac{2\left(\bar{Y}+B_\phi P\right)}{P}Cov\left(p,\bar{y}\right) + V\left(\bar{y}\right) \qquad (2.4)$$

Simplifying (2.4), we get MSE of $t_1$ as

$$MSE\left(t_1\right) \cong \left(\frac{1-f}{n}\right)\left[R_1^2S_\phi^2 + S_y^2\left(1-\rho_{pb}^2\right)\right] \qquad (2.5)$$

Several authors have used prior value of certain population parameters (s) to find more precise estimates. Searls (1964) used Coefficient of Variation (CV) of study character at estimation stage. In practice this CV is seldom known. Motivated by Searls (1964) work, Sen (1978), Sisodiya and Dwivedi (1981), and Upadhyaya and Singh (1984) used the known CV of the auxiliary character for estimating population mean of a study character in ratio method of estimation. The use of prior value of Coefficient of Kurtosis in estimating the population variance of study character y was first made by Singh et. al. (1973). Later, used by and Searls and Intarapanich (1990), Upadhyaya and Singh (1999), Singh (2003) and Singh et. al. (2004) in the estimation of population mean of study character. Recently Singh and Tailor (2003) proposed a modified ratio estimator by using the known value of correlation coefficient.

In next section, we propose some ratio estimators for estimating the population mean of the variable under study using known parameters of the attribute $\phi$ such as coefficient of variation $C_p$, Kurtosis $(\beta_2(\phi))$ and point biserial correlation coefficient $\rho_{pb}$.

## 3. Suggested Estimators

We suggest following estimator –

$$t = \frac{\bar{y} + b_\phi(P-p)}{(m_1 p + m_2)}(m_1 P + m_2) \tag{3.1}$$

Where $m_1(\neq 0)$, $m_2$ are either real number or the functions of the known parameters of the attribute such as $C_p$, $(\beta_2(\phi))$ and $\rho_{pb}$.

The following scheme presents some of the important estimators of the population mean, which can be obtained by suitable choice of constants $m_1$ and $m_2$:

| Estimator | Values of | |
|---|---|---|
| | $m_1$ | $m_2$ |
| $t_1 = \frac{\bar{y} + b_\phi(P-p)}{p}P$ | 1 | 1 |
| $t_2 = \frac{\bar{y} + b_\phi(P-p)}{(p + \beta_2(\phi))}[P + \beta_2(\phi)]$ | 1 | $\beta_2(\phi)$ |
| $t_3 = \frac{\bar{y} + b_\phi(P-p)}{(p + C_p)}[P + C_p]$ | 1 | $C_p$ |
| $t_4 = \frac{\bar{y} + b_\phi(P-p)}{(p + \rho_{pb}))}[P + \rho_{pb}]$ | 1 | $\rho_{pb}$ |
| $t_5 = \frac{\bar{y} + b_\phi(P-p)}{(p\beta_2(\phi) + C_p)}[P\beta_2(\phi) + C_p]$ | $\beta_2(\phi)$ | $C_p$ |
| $t_6 = \frac{\bar{y} + b_\phi(P-p)}{(pC_p + \beta_2(\phi))}[PC_p + \beta_2(\phi)]$ | $C_p$ | $\beta_2(\phi)$ |
| $t_7 = \frac{\bar{y} + b_\phi(P-p)}{(pC_p + \rho_{pb})}[PC_p + \rho_{pb}]$ | $C_p$ | $\rho_{pb}$ |
| $t_8 = \frac{\bar{y} + b_\phi(P-p)}{(p\rho_{pb} + C_p)}[P\rho_{pb} + C_p]$ | $\rho_{pb}$ | $C_p$ |
| $t_9 = \frac{\bar{y} + b_\phi(P-p)}{(p\beta_2(\phi) + \rho_{pb})}[P\beta_2(\phi) + \rho_{pb}]$ | $\beta_2(\phi)$ | $\rho_{pb}$ |
| $t_{10} = \frac{\bar{y} + b_\phi(P-p)}{(p\rho_{pb} + \beta_2(\phi))}[P\rho_{pb} + \beta_2(\phi)]$ | $\rho_{pb}$ | $\beta_2(\phi)$ |

Following the approach of section 2, we obtain the MSE expression for these proposed estimators as –

$$MSE(t_i) \cong \left(\frac{1-f}{n}\right)\left[R_i S_\phi^2 + S_y^2(1-\rho_{pb}^2)\right], \qquad (i = 1,2,3,....,10) \qquad (3.2)$$

Where $R_1 = \frac{\overline{Y}}{P}$, $R_2 = \frac{\overline{Y}}{P+\beta_2(\phi)}$, $R_3 = \frac{\overline{Y}}{P+C_p}$, $R_4 = \frac{\overline{Y}}{P+\rho_{pb}}$,

$R_5 = \frac{\overline{Y}\beta_2(\phi)}{P\beta_2(\phi)+C_p}$, $R_6 = \frac{\overline{Y}C_p}{PC_p+\beta_2(\phi)}$, $R_7 = \frac{\overline{Y}C_p}{PC_p+\rho_{pb}}$, $R_8 = \frac{\overline{Y}\rho_{pb}}{P\rho_{pb}+C_p}$,

$R_9 = \frac{\overline{Y}\beta_2(\phi)}{P\beta_2(\phi)+\rho_{pb}}$, $R_{10} = \frac{\overline{Y}\rho_{pb}}{P\rho_{pb}+\beta_2(\phi)}$.

## 4. Efficiency comparisons

It is well known that under simple random sampling without replacement (SRSWOR) the variance of the sample mean is

$$V(\overline{y}) = \left(\frac{1-f}{n}\right)S_y^2 \qquad (4.1)$$

From (4.1) and (3.2), we have

$$V(\overline{y}) - MSE(t_i) \geq 0, \qquad i = 1,2,.....,10$$

$$\Rightarrow \rho_{pb}^2 > \frac{S_\phi^2}{S_y^2}R_i^2 \qquad (4.2)$$

When this condition is satisfied, proposed estimators are more efficient than the sample mean.

Now, we compare the MSE of the proposed estimators with the MSE of Naik and Gupta (1996) estimator $t_{NG}$. From (3.2) and (1.1) we have

$$MSE(t_{NG}) - MSE(t_i) \geq 0, \qquad (i = 1,2,.....,10)$$

$$\Rightarrow \rho_{pb}^2 \geq \frac{S_\phi^2}{S_y^2}\left[R_i^2 - R_\phi^2 + 2R_\phi K_{yp}\right] \qquad (4.3)$$

where $K_{yp} = \rho_{yp}\frac{C_y}{C_p}$.

## 5. Empirical Study

The data for the empirical study is taken from natural population data set considered by Sukhatme and Sukhatme (1970):

$y$ = Number of villages in the circles and
$\phi$ = A circle consisting more than five villages

N = 89, $\bar{Y}$ = 3.36, P = 0.1236, $\rho_{pb}$ = 0.766, $C_y$ = 0.604, $C_p$ = 2.19, $\beta_2(\phi)$ = 6.23181.

In table 5.1 percent relative efficiencies (PRE) of various estimators are computed with respect to $\bar{y}$.

**Table 5.1: PRE of different estimators of $\bar{Y}$ with respect to $\bar{y}$**

| Estimator | PRE (., $\bar{y}$) |
|---|---|
| $\bar{y}$ | 100 |
| $t_{NG}$ | 11.61 |
| $t_1$ | 7.36 |
| $t_2$ | 236.55 |
| $t_3$ | 227.69 |
| $t_4$ | 208.09 |
| $t_5$ | 185.42 |
| $t_6$ | 230.72 |
| $t_7$ | 185.27 |
| $t_8$ | 230.77 |
| $t_9$ | 152.37 |
| $t_{10}$ | 237.81 |

From table 5.1, we observe that the proposed estimators $t_i(i = 2,.....,10)$ which uses some known values of population proportion performs better than the usual sample mean $\bar{y}$ and Naik and Gupta (1996) estimator $t_{NG}$.

## Conclusion

We have suggested some ratio estimators for estimating $\bar{Y}$ which uses some known value of population proportion. For practical purposes the choice of the estimator depends upon the availability of the population parameters.